\input amstex
\documentstyle{amsppt}
%
\catcode`@=11
\redefine\output@{%
  \def\break{\penalty-\@M}\let\par\endgraf
  \ifodd\pageno\global\hoffset=105pt\else\global\hoffset=8pt\fi  
  \shipout\vbox{%
    \ifplain@
      \let\makeheadline\relax \let\makefootline\relax
    \else
      \iffirstpage@ \global\firstpage@false
        \let\rightheadline\frheadline
        \let\leftheadline\flheadline
      \else
        \ifrunheads@ 
        \else \let\makeheadline\relax
        \fi
      \fi
    \fi
    \makeheadline \pagebody \makefootline}%
  \advancepageno \ifnum\outputpenalty>-\@MM\else\dosupereject\fi
}
\def\Beta{\mathchar"0\hexnumber@\rmfam 42}
\catcode`\@=\active
\nopagenumbers
\def\negskp{\hskip -2pt}

\accentedsymbol\hatgamma{\kern 2pt\hat{\kern -2pt\gamma}}
\accentedsymbol\checkgamma{\kern 2.5pt\check{\kern -2.5pt\gamma}}
\def\blue#1{#1}

\catcode`#=11\def\diez{#}\catcode`#=6
\catcode`&=11\catcode`&=4
\catcode`_=11\def\podcherkivanie{_}\catcode`_=8
\def\mycite#1{\cite{\blue{#1}}\immediate\special{ps:
     ShrHPSdict begin /ShrBORDERthickness 0 def}}
\def\myciterange#1#2#3#4{\cite{\blue{#2#3#4}}\immediate\special{ps:
     ShrHPSdict begin /ShrBORDERthickness 0 def}}
\def\mytag#1{%
    \tag#1}
\def\mythetag#1{\thetag{\blue{#1}}\immediate\special{ps:
     ShrHPSdict begin /ShrBORDERthickness 0 def}}
\def\myrefno#1{\no#1}
\def\myhref#1#2{\blue{#2}\immediate\special{ps:
     ShrHPSdict begin /ShrBORDERthickness 0 def}}
\def\myEarXivlink{\myhref{http://arXiv.org}{http:/\negskp/arXiv.org}}

\def\mytheorem#1{\csname proclaim\endcsname{Theorem #1}}
\def\mytheoremwithtitle#1#2{\csname proclaim\endcsname{Theorem #1#2}}
\def\mythetheorem#1{\blue{#1}\immediate\special{ps:
     ShrHPSdict begin /ShrBORDERthickness 0 def}}
\def\mylemma#1{\csname proclaim\endcsname{Lemma #1}}
\def\mylemmawithtitle#1#2{\csname proclaim\endcsname{Lemma #1#2}}
\def\mythelemma#1{\blue{#1}\immediate\special{ps:
     ShrHPSdict begin /ShrBORDERthickness 0 def}}
\def\mycorollary#1{\csname proclaim\endcsname{Corollary #1}}

\def\myconjecture#1{\csname proclaim\endcsname{Conjecture #1}}
\def\myconjecturewithtitle#1#2{\csname proclaim\endcsname{Conjecture #1#2}}
\def\mytheconjecture#1{\blue{#1}\immediate\special{ps:
     ShrHPSdict begin /ShrBORDERthickness 0 def}}

\pagewidth{360pt}
\pageheight{606pt}
\topmatter
\title
A note on the second
cuboid\\ conjecture. Part~\uppercase\expandafter{\romannumeral 1}.
\endtitle
\author
Ruslan Sharipov
\endauthor
\address Bashkir State University, 32 Zaki Validi street, 450074 Ufa, Russia
\endaddress
\email\myhref{mailto:r-sharipov\@mail.ru}{r-sharipov\@mail.ru}
\endemail
\abstract
    The problem of finding perfect Euler cuboids or proving their non-existence
is an old unsolved problem in mathematics. The second cuboid conjecture is one of 
the three propositions suggested as intermediate stages in proving the 
non-existence of perfect Euler cuboids. It is associated with a certain polynomial 
Diophantine equation of the order 10. In this paper a structural theorem for the 
solutions of this Diophantine equation is proved and some examples of its
application are considered. 
\endabstract
\subjclassyear{2000}
\subjclass 11D41, 11D72, 12E05\endsubjclass
\endtopmatter
\TagsOnRight
\document

\head
1. Introduction.
\endhead
     Let's denote through $Q_{p\kern 0.6pt q}(t)$ the following polynomial 
of the order $10$ depending on two integer parameters $p$ and $q$:
$$
\hskip -2em
\gathered
Q_{p\kern 0.6pt q}(t)=t^{10}+(2\,q^2+p^2)\,(3\,q^2-2\,p^2)\,t^8
+(q^8+10\,p^2\,q^6+\\
+\,4\,p^4\,q^4-14\,p^6\,q^2+p^8)\,t^6
-p^2\,q^2\,(p^8+10\,q^2\,p^6+4\,q^4\,p^4-\\
-\,14\,q^6\,p^2+q^8)\,t^4-p^6\,q^6\,(2\,p^2+q^2)
\,(3\,p^2-2\,q^2)\,t^2-q^{10}\,p^{10}.
\endgathered
\mytag{1.1}
$$
\myconjecturewithtitle{1.1}{ (second cuboid conjecture)} For any positive 
coprime integers $p\neq q$ the polynomial $Q_{p\kern 0.6pt q}(t)$ in 
\mythetag{1.1} is irreducible in the ring $\Bbb Z[t]$. 
\endproclaim
     The second cuboid conjecture~\mytheconjecture{1.1} and the polynomial
$Q_{p\kern 0.6pt q}(t)$ in it were introduced in \mycite{1}. They are 
associated with the problem of constructing a perfect Euler cuboid
(see \mycite{2} and \myciterange{3}{3}{--}{37} for more details). Let's
write the equation
$$
\hskip -2em
Q_{p\kern 0.6pt q}(t)=0.
\mytag{1.2}
$$
The equation \mythetag{1.2} can be understood as a Diophantine equation of 
the order 10 with two integer parameters $p$ and $q$. The second cuboid 
conjecture~\mytheconjecture{1.1} implies the following theorem.
\mytheorem{1.1} For any positive coprime integers $p\neq q$ the polynomial 
Diophantine equation \mythetag{1.2} has no integer solutions. 
\endproclaim
\noindent Note that a similar theorem associated with the first cuboid 
conjecture was formulated and proved in \mycite{38}.\par
     The theorem~\mythetheorem{1.1} is a weaker proposition than the 
conjecture~\mytheconjecture{1.1} itself. However, even proving this proposition
\pagebreak in the case of the second cuboid conjecture is rather difficult. 
Below in section~4 we formulate and prove a structural theorem for the 
solutions of the Diophantine equation \mythetag{1.2}, if any, and in section~5 we
use it in order to prove the theorem~\mythetheorem{1.1} for some particular values 
$p$ and $q$. 
\head
2. The inversion symmetry and parity.
\endhead
     The polynomial $Q_{p\kern 0.6pt q}(t)$ in \mythetag{1.1} possesses some 
special properties. They are expressed by the following formulas which can be 
verified by direct calculations:
$$
\xalignat 2
&\hskip -2em
Q_{p\kern 0.6pt q}(t)=-\frac{Q_{qp}(p^2\,q^2/t)\,t^{10}}{p^{10}\,q^{10}},
&&Q_{p\kern 0.6pt q}(t)=Q_{p\kern 0.6pt q}(-t),\\
\vspace{-1.3ex}
\mytag{2.1}\\
\vspace{-1.3ex}
&\hskip -2em
Q_{qp}(t)=-\frac{Q_{p\kern 0.6pt q}(p^2\,q^2/t)\,t^{10}}{p^{10}\,q^{10}},
&&Q_{qp}(t)=Q_{qp}(-t).
\endxalignat
$$
Note that in \mythetag{2.1} we have two polynomials $Q_{p\kern 0.6pt q}(t)$
and $Q_{qp}(t)$. The polynomial $Q_{qp}(t)$ is produced from \mythetag{1.1}
by exchanging parameters $p$ and $q$:
$$
\hskip -2em
\gathered
Q_{qp}(t)=t^{10}+(2\,p^2+q^2)\,(3\,p^2-2\,q^2)\,t^8
+(p^8+10\,q^2\,p^6+\\
+\,4\,q^4\,p^4-14\,q^6\,p^2+q^8)\,t^6
-q^2\,p^2\,(q^8+10\,p^2\,q^6+4\,p^4\,q^4+\\
-\,14\,p^6\,q^2+p^8)\,t^4-q^6\,p^6\,(2\,q^2+p^2)
\,(3\,q^2-2\,p^2)\,t^2-p^{10}\,q^{10}.
\endgathered
\mytag{2.2}
$$\par
     Two of the four symmetries in \mythetag{2.1} contain the inversion of $t$.
For this reason they are called inversion symmetries. The other two symmetries
in \mythetag{2.1} mean that the polynomials \mythetag{1.1} and \mythetag{2.2}
are even with respect to their argument $t$.
\head
3. Some prerequisites.
\endhead
     Assume that the polynomial $Q_{p\kern 0.6pt q}(t)$ has an integer root 
$t=A_0$. Since $p\neq 0$ and $q\neq 0$, we have $A_0\neq 0$. Then due to the 
inversion symmetries in \mythetag{2.1} the polynomial $Q_{qp}(t)$ has an integer 
root $t=B_0$, where
$$
\hskip -2em
B_0=\frac{p^2\,q^2}{A_0}. 
\mytag{3.1}
$$
Since $p\neq 0$ and $q\neq 0$, from \mythetag{3.1} we derive $B_0\neq 0$. 
Applying the parity symmetry from \mythetag{2.1}, we conclude that the polynomial
$Q_{p\kern 0.6pt q}(t)$ has the other integer root\linebreak $t=-A_0$, while 
$Q_{qp}(t)$ has the other integer root $t=-B_0$. As a result the polynomials 
$Q_{p\kern 0.6pt q}(t)$ and $Q_{qp}(t)$ split into factors
$$
\xalignat 2
&\hskip -2em
Q_{p\kern 0.6pt q}(t)=(t^2-A_0^2)\ C_8(t),
&&Q_{qp}(t)=(t^2-B_0^2)\ D_8(t)\quad
\mytag{3.2}
\endxalignat
$$
with $A_0>0$ and $B_0>0$. Here $C_8(t)$ and $D_8(t)$ are eighth order 
polynomials complementary to $t^2-A_0^2$ and $t^2-B_0^2$. Applying 
\mythetag{2.1} to \mythetag{3.2} we derive
$$
\xalignat 2
&\hskip -2em
C_8(t)=C_8(-t),
&&D_8(t)=D_8(-t).\quad
\mytag{3.3}
\endxalignat
$$
Due to \mythetag{3.3} the polynomials $C_8(t)$ and $D_8(t)$ are given by the
formulas
$$
\hskip -2em
\aligned
&C_8(t)=t^8+C_6\,t^6+C_4\,t^4+C_2\,t^2+C_0,\\
&D_8(t)=t^8+D_6\,t^6+D_4\,t^4+D_2\,t^2+D_0.
\endaligned
\mytag{3.4}
$$
The coefficients of the polynomials \mythetag{3.4} are integer numbers.\par
     Now let's apply the inversion symmetries from \mythetag{2.1} to
\mythetag{3.2}. As a result we get
$$
\xalignat 2
&\hskip -2em
C_8(t)=\frac{D_8(p^2\,q^2/t)\,t^8}{p^6\,q^6\,A_0^2},
&&D_8(t)=\frac{C_8(p^2\,q^2/t)\,t^8}{p^6\,q^6\,B_0^2}.
\mytag{3.5}
\endxalignat
$$
Applying the symmetries \mythetag{3.5} to \mythetag{3.4}, we derive a series 
of relationships for the coefficients of the polynomials $C_8(t)$ and
$D_8(t)$:
$$
\xalignat 2
&\hskip -2em
C_0\,A_0^2=p^{10}\,q^{10},
&&C_2\,A_0^2=p^6\,q^6\,D_6,\\
&\hskip -2em
C_4\,A_0^2=p^2\,q^2\,D_4,
&&C_6\,A_0^2\,p^2\,q^2=D_2,\\
&\hskip -2em
A_0^2\,p^6\,q^6=D_0,
&&D_0\,B_0^2=p^{10}\,q^{10},
\mytag
{3.6}\\
&\hskip -2em
D_2\,B_0^2=p^6\,q^6\,C_6,
&&D_4\,B_0^2=p^2\,q^2\,C_4,\\
&\hskip -2em
D_6\,B_0^2\,p^2\,q^2=C_2,
&&B_0^2\,p^6\,q^6=C_0.
\endxalignat
$$
The equations \mythetag{3.6} are excessive. Due to \mythetag{3.1} some of them 
are equivalent to some others. For this reason we can eliminate excessive 
variables:
$$
\xalignat 2
&\hskip -2em
D_0=p^6\,q^6\,A_0^2,
&&D_2=C_6\,p^2\,q^2\,A_0^2,\\
\vspace{-1.5ex}
\mytag{3.7}\\
\vspace{-1.5ex}
&\hskip -2em
C_0=p^6\,q^6\,B_0^2,
&&C_2=D_6\,p^2\,q^2\,B_0^2.
\endxalignat
$$
Substituting \mythetag{3.7} into the formulas \mythetag{3.4} for $C_8(t)$
and $D_8(t)$, we get 
$$
\hskip -2em
\aligned
&C_8(t)=t^8+C_6\,t^6+C_4\,t^4+D_6\,p^2\,q^2\,B_0^2\,t^2+p^6\,q^6\,B_0^2,\\
&D_8(t)=t^8+D_6\,t^6+D_4\,t^4+C_6\,p^2\,q^2\,A_0^2\,t^2+p^6\,q^6\,A_0^2.
\endaligned
\mytag{3.8}
$$\par
     Unlike $C_0$, $D_0$, $C_2$, and $D_2$ in \mythetag{3.7}, the coefficients
$C_4$ and $D_4$ in \mythetag{3.8} are not expressed through other coefficients. 
However, they are not independent. They are related with each other by means
of the equation 
$$
\hskip -2em
A_0\,C_4=B_0\,D_4.
\mytag{3.9}
$$
The equation \mythetag{3.9} is derived from \mythetag{3.6} by means of the
formula \mythetag{3.1}.\par
     Having derived the formulas \mythetag{3.8}, we substitute them back into
the relationships \mythetag{3.2}. As a result we derive the following formulas:
$$
\gather
\hskip -2em
\gathered
Q_{p\kern 0.6pt q}(t)=t^{10}+(C_6-A_0^2)\,t^8+(C_4-A_0^2\,C_6)\,t^6
+(D_6\,p^2\,q^2\,B_0^2\,-\\
-\,A_0^2\,C_4)\,t^4+q^2\,p^2\,B_0^2\,(p^4\,q^4-A_0^2\,D_6)\,t^2
-A_0^2\,p^6\,q^6\,B_0^2,
\endgathered
\mytag{3.10}\\
\vspace{1ex}
\hskip -2em
\gathered
Q_{qp}=t^{10}+(D_6-B_0^2)\,t^8+(D_4-B_0^2\,D_6)\,t^6+(C_6\,p^2\,q^2\,A_0^2\,-\\
-\,B_0^2\,D_4)\,t^4+p^2\,q^2\,A_0^2\,(p^4\,q^4-B_0^2\,C_6)\,t^2
-A_0^2\,p^6\,q^6\,B_0^2.
\endgathered
\mytag{3.11}
\endgather
$$
Comparing the formula \mythetag{3.10} with \mythetag{1.1} and comparing 
the formula \mythetag{3.11} with \mythetag{2.2}, we derive ten equations
for the coefficients of the polynomials \mythetag{3.8}. Two of them are
equivalent to the equation \mythetag{3.1} written as
$$
\hskip -2em
A_0\,B_0=p^2\,q^2. 
\mytag{3.12}
$$ 
The other eight equations are written as follows:
$$
\gather
\hskip -2em
\gathered
C_6-A_0^2=(2\,q^2+p^2)\,(3\,q^2-2\,p^2),\\
D_6-B_0^2=(2\,p^2+q^2)\,(3\,p^2-2\,q^2),
\endgathered
\mytag{3.13}\\
\vspace{1ex}
\hskip -2em
\gathered
C_4-A_0^2\,C_6=p^8-14\,p^6\,q^2+4\,p^4\,q^4+10\,p^2\,q^6+q^8,\\
D_4-B_0^2\,D_6=q^8-14\,q^6\,p^2+4\,q^4\,p^4+10\,q^2\,p^6+p^8,
\endgathered
\mytag{3.14}\\
\vspace{1ex}
\hskip -2em
\gathered
A_0^2\,C_4-D_6\,p^2\,q^2\,B_0^2=p^2\,q^2\,(q^8-14\,q^6\,p^2+4\,q^4\,p^4
+10\,q^2\,p^6+p^8),\\
B_0^2\,D_4-C_6\,p^2\,q^2\,A_0^2=p^2\,q^2\,(p^8-14\,p^6\,q^2+4\,p^4\,q^4
+10\,p^2\,q^6+q^8),
\endgathered\quad
\mytag{3.15}\\
\vspace{1ex}
\hskip -2em
\gathered
B_0^2\,(A_0^2\,D_6-p^4\,q^4)=p^4\,q^4\,(2\,p^2+q^2)\,(3\,p^2-2\,q^2),\\
A_0^2\,(B_0^2\,C_6-p^4\,q^4)=p^4\,q^4\,(2\,q^2+p^2)\,(3\,q^2-2\,p^2).
\endgathered
\mytag{3.16}
\endgather
$$
The equations \mythetag{3.12}, \mythetag{3.13}, \mythetag{3.14}, 
\mythetag{3.15}, \mythetag{3.16} are excessive. Indeed, the equations
\mythetag{3.16} follow from \mythetag{3.12} and \mythetag{3.13}. 
Similarly, the equations \mythetag{3.15} can be derived from 
\mythetag{3.14} with the use of \mythetag{3.9} and \mythetag{3.12}. 
As for the equations \mythetag{3.13} and \mythetag{3.14}, when
complemented with the equations \mythetag{3.9} and \mythetag{3.12},
they constitute a system of Diophantine equations with respect to
$C_4$, $D_4$, $C_6$, $D_6$, $A_0$ and $B_0$. The results of the above 
calculations are summarized in the following lemma.
\mylemma{3.1} For $p\neq 0$ and $q\neq 0$ the polynomial 
$Q_{p\kern 0.6pt q}(t)$ in \mythetag{1.1} has integer roots if and only 
if the system of Diophantine equations \mythetag{3.9}, \mythetag{3.12},
\mythetag{3.13}, and \mythetag{3.14} is solvable with respect to
the integer variables $C_4$, $D_4$, $C_6$, $D_6$, $A_0>0$ and $B_0>0$.
\endproclaim 
\head
4. The structural theorem.
\endhead
    Below we continue studying the equations \mythetag{3.9}, 
\mythetag{3.12}, \mythetag{3.13}, \mythetag{3.14} implicitly assuming 
$p\neq q$ to be two positive coprime integer numbers. Let $p_1,\,\ldots,\,
p_m$ be the prime factors of $p$ and let $q_1,\,\ldots,\,q_n$ be
the prime factors of $q$:
$$
\xalignat 2
&\hskip -2em
p=p_1^{\kern 0.4pt\alpha_1}\cdot\ldots\cdot p_m^{\kern 0.4pt\alpha_m},
&&q=q_1^{\kern 0.4pt\beta_1}\cdot\ldots\cdot q_n^{\kern 0.4pt\beta_n}.
\mytag{4.1}
\endxalignat
$$ 
Usually the multiplicities $\alpha_1,\,\ldots,\,\alpha_m$ and
$\beta_1,\,\ldots,\,\beta_n$ in \mythetag{4.1} are positive numbers. 
However, in order to cover two special cases $p=1$ and $q=1$ we assume 
them to be non-negative numbers. Since $p$ and $q$ are assumed coprime,
i\.\,e\. 
$$
\hskip -2em
\gcd(p,q)=1,
\mytag{4.2}
$$
the prime factors $p_1,\,\ldots,\,p_m$ and $q_1,\,\ldots,\,q_n$ in
\mythetag{4.1} are distinct, i\.\,e\. $p_i\neq q_j$. 
\mylemma{4.1} For any solution of the Diophantine equations \mythetag{3.9}, 
\mythetag{3.12}, \mythetag{3.13}, and \mythetag{3.14} with $A_0>0$ and $B_0>0$ 
if\/ $A_0\neq 1$, each prime factor\/ $r$ of\/ $A_0$ is a prime factor of\/ $p$ 
or a prime factor of\/ $q$, i\.\,e\. $r=p_i$ or\/ $r=q_j$. 
\endproclaim
\mylemma{4.2} For any solution of the Diophantine equations \mythetag{3.9}, 
\mythetag{3.12}, \mythetag{3.13}, and \mythetag{3.14} with $A_0>0$ and $B_0>0$
if\/ $B_0\neq 1$, each prime factor\/ $r$ of\/ $B_0$ is a prime factor of\/ $p$ 
or a prime factor of\/ $q$, i\.\,e\. $r=p_i$ or\/ $r=q_j$. 
\endproclaim
    The lemmas~\mythelemma{4.1} and \mythelemma{4.2} are immediate from
\mythetag{4.2} and \mythetag{3.12}. \pagebreak Due to the 
lemmas~\mythelemma{4.1} and \mythelemma{4.2} we can write the following 
expansions for $A_0$ and $B_0$:
$$
\hskip -2em
\aligned
&A_0=p_1^{\kern 0.7pt\mu_1}\cdot\ldots\cdot p_m^{\kern 0.7pt\mu_m}
\,q_1^{\nu_1}\cdot\ldots\cdot q_n^{\nu_n},\\
&B_0=p_1^{\kern 0.7pt\eta_1}\cdot\ldots\cdot p_m^{\kern 0.7pt\eta_m}
\,q_1^{\tau_1}\cdot\ldots\cdot q_n^{\tau_n}.
\endaligned
\mytag{4.3}
$$
The multiplicities $\mu_1,\,\ldots,\,\mu_m$, $\nu_1,\,\ldots,\,\nu_n$, 
$\eta_1,\,\ldots,\,\eta_m$, and $\tau_1,\,\ldots,\,\tau_n$ in the
expansions \mythetag{4.3} obey the following relationships: 
$$
\xalignat 2
&\hskip -2em
\mu_i+\eta_i=2\,\alpha_i,
&&\nu_i+\tau_i=2\,\beta_i.
\mytag{4.4}
\endxalignat
$$
The formulas \mythetag{4.4} are easily derived by substituting
the expansions \mythetag{4.1} and \mythetag{4.3} into the equation 
\mythetag{3.12}.\par
     Now let's return back to the equations \mythetag{3.13} and
\mythetag{3.14}. It is easy to see that the equations \mythetag{3.13}
can be explicitly resolved with respect to $C_6$ and $D_6$:
$$
\hskip -2em
\gathered
C_6=A_0^2+(2\,q^2+p^2)\,(3\,q^2-2\,p^2),\\
D_6=B_0^2+(2\,p^2+q^2)\,(3\,p^2-2\,q^2).
\endgathered
\mytag{4.5}
$$
Upon substituting \mythetag{4.5} into the equations \mythetag{3.14} we 
can explicitly resolve the equations \mythetag{3.14} with respect to the 
variables $C_4$ and $D_4$:
$$
\align
&\hskip -2em
\aligned
C_4&=A_0^4+(2\,q^2+p^2)\,(3\,q^2-2\,p^2)\,A_0^2\,+\\
&+p^8-14\,p^6\,q^2+4\,p^4\,q^4+10\,p^2\,q^6+q^8,
\endaligned
\mytag{4.6}\\
\vspace{1ex}
&\hskip -2em
\aligned
D_4&=B_0^4+(2\,p^2+q^2)\,(3\,p^2-2\,q^2)\,B_0^2\,+\\
&+q^8-14\,q^6\,p^2+4\,q^4\,p^4+10\,q^2\,p^6+p^8.
\endaligned
\mytag{4.7}
\endalign
$$
And finally, we can substitute \mythetag{4.6} and \mythetag{4.7} into
the equation \mythetag{3.9}. As a result we derive the following equation
for the variables $A_0$ and $B_0$: 
$$
\hskip -2em
\aligned
A_0\,&\bigl(A_0^4+(2\,q^2+p^2)\,(3\,q^2-2\,p^2)\,A_0^2\,+\\
&+p^8-14\,p^6\,q^2+4\,p^4\,q^4+10\,p^2\,q^6+q^8\bigr)=\\
=\,\,&B_0\,\bigl(B_0^4+(2\,p^2+q^2)\,(3\,p^2-2\,q^2)\,B_0^2\,+\\
&+q^8-14\,q^6\,p^2+4\,q^4\,p^4+10\,q^2\,p^6+p^8\bigr).
\endaligned
\mytag{4.8}
$$
Summarizing these calculations, we can formulate the following lemma.
\mylemma{4.3} The system of four Diophantine equations \mythetag{3.9}, 
\mythetag{3.12}, \mythetag{3.13}, \mythetag{3.14} is equivalent to
the system of two Diophantine equations \mythetag{3.12} and \mythetag{4.8}.
\endproclaim
\mylemma{4.4} For any solution of the Diophantine equations \mythetag{3.9}, 
\mythetag{3.12}, \mythetag{3.13}, and \mythetag{3.14} with $A_0>0$ and 
$B_0>0$ if\/ $\mu_i>0$ and $\eta_i>0$ in \mythetag{4.3}, then $\mu_i=\eta_i
=\alpha_i$.
\endproclaim
\demo{Proof} The proof is by contradiction. Assume that $\mu_i>0$, 
$\eta_i>0$, and $\mu_i\neq\eta_i$. Then from \mythetag{4.6} and 
\mythetag{4.7}, applying \mythetag{4.1}, \mythetag{4.2}, \mythetag{4.3}, 
and \mythetag{4.4}, we derive 
$$
\xalignat 2
&\hskip -2em
C_4\equiv q^8\ (\kern -0.6em\mod p_i),
&&D_4\equiv q^8\ (\kern -0.6em\mod p_i).
\mytag{4.9}
\endxalignat
$$
Moreover, the formula \mythetag{4.3} yields the relationships
$$
\xalignat 2
&\hskip -2em
A_0=A'_0\,p^{\kern 0.7pt\mu_i}_i,
&&B_0=B'_0\,p^{\kern 0.7pt\eta_i}_i,
\mytag{4.10}
\endxalignat
$$
where $A'_0\not\equiv 0\,(\kern -0.6em\mod p_i)$ and $B'_0\not\equiv 0\,
(\kern -0.6em\mod p_i)$. Our assumption $\mu_i\neq\eta_i$ means that 
$\mu_i>\eta_i$ or $\mu_i<\eta_i$. If\/ $\mu_i>\eta_i$, then substituting 
\mythetag{4.10} into \mythetag{3.9}, we derive 
$$
\hskip -2em
A_0'\,C_4\,p^{\kern 0.7pt\mu_i-\eta_i}_i=B'_0\,D_4.
\mytag{4.11}
$$
Due to \mythetag{4.10}, \mythetag{4.9} and \mythetag{4.2} the left hand side 
of \mythetag{4.11} is zero modulo $p_i$, while the right hand side of 
\mythetag{4.11} is nonzero modulo $p_i$, which is contradictory.\par
     Similarly, if\/ $\mu_i<\eta_i$, substituting \mythetag{4.10} into 
\mythetag{3.9}, we derive
$$
\hskip -2em
A_0'\,C_4=B'_0\,D_4\,p^{\kern 0.7pt\eta_i-\mu_i}_i.
\mytag{4.12}
$$
In this case the left hand side of \mythetag{4.12} is nonzero modulo $p_i$, 
while the right hand side of \mythetag{4.12} is zero modulo $p_i$, which is 
also contradictory.\par
     The contradictions obtained prove that $\mu_i=\eta_i$. The equalities
$\mu_i=\alpha_i$ and $\eta_i=\alpha_i$ are immediate from $\mu_i=\eta_i$ due
to \mythetag{4.4}. The lemma~\mythelemma{4.4} is proved.
\qed\enddemo
\mylemma{4.5} For any solution of the Diophantine equations \mythetag{3.9}, 
\mythetag{3.12}, \mythetag{3.13}, and \mythetag{3.14} with $A_0>0$ and 
$B_0>0$ if\/ $\nu_i>0$ and $\tau_i>0$ in \mythetag{4.3}, then $\nu_i=\tau_i
=\beta_i$.
\endproclaim
     The lemma~\mythelemma{4.5} is analogous to the lemma~\mythelemma{4.4}.
Its proof is similar to the above proof of the lemma~\mythelemma{4.4}.\par
     Now, relying on the lemmas~\mythelemma{4.4} and \mythelemma{4.5}, we 
define the following integer numbers:
$$
\xalignat 3
&\hskip -2em
a_p\,=\!\!\prod_{\eta_i=0}\!p^{\kern 0.4pt\alpha_i}_i,
&&b_p\,=\!\!\prod_{\mu_i=0}\!p^{\kern 0.4pt\alpha_i}_i,
&&c_p\,=\!\!\prod\Sb\mu_i>0\\ \eta_i>0\endSb\!p^{\kern 0.4pt\alpha_i}_i,
\quad
\mytag{4.13}\\
\vspace{2ex}
&\hskip -2em
a_q\,=\!\!\prod_{\tau_i=0}\!q^{\kern 0.4pt\beta_i}_i,
&&b_q\,=\!\!\prod_{\nu_i=0}\!q^{\kern 0.4pt\beta_i}_i,
&&c_q\,=\!\!\prod\Sb\nu_i>0\\ \tau_i>0\endSb\!q^{\kern 0.4pt\beta_i}_i.
\quad
\mytag{4.14}
\endxalignat
$$
Note that the numbers \mythetag{4.13} and \mythetag{4.14} are pairwise mutually
coprime, i\.\,e\.
$$
\xalignat 3
&\hskip -2em
\gcd(a_p,b_p)=1, &&\gcd(a_p,c_p)=1, &&\gcd(a_p,a_q)=1,\\
&\hskip -2em
\gcd(a_p,b_q)=1, &&\gcd(a_p,c_q)=1, &&\gcd(b_p,c_p)=1,\\
&\hskip -2em
\gcd(b_p,a_q)=1, &&\gcd(b_p,b_q)=1, &&\gcd(b_p,c_q)=1,
\quad\mytag{4.15}\\
&\hskip -2em
\gcd(c_p,a_q)=1, &&\gcd(c_p,b_q)=1, &&\gcd(c_p,c_q)=1,\\
&\hskip -2em
\gcd(a_q,b_q)=1, &&\gcd(a_q,c_q)=1, &&\gcd(b_q,c_q)=1.
\endxalignat
$$
If\/ $\mu_i=0$, then $\eta_i=2\,\alpha_i$, and if\/ $\eta_i=0$, then  
$\mu_i=2\,\alpha_i$. Similarly, if\/ $\nu_i=0$, then $\tau_i=2\,\beta_i$
and if\/ $\tau_i=0$, then $\nu_i=2\,\beta_i$. These implications are 
derived from \mythetag{4.4}. The lemmas~\mythelemma{4.4} and \mythelemma{4.5} 
say that if $\mu_i>0$ and $\eta_i>0$, then $\mu_i=\eta_i
=\alpha_i$, and if $\nu_i>0$ and $\tau_i>0$, then $\nu_i=\tau_i=\beta_i$. 
As a result from \mythetag{4.13} and \mythetag{4.14} we derive
$$
\pagebreak
\xalignat 2
\kern -2em
A_0&=a_p^2\,c_p\,a_q^2\,c_q,
&B_0&=b_p^2\,c_p\,b_q^2\,c_q,
\qquad\\
\vspace{-1.5ex}
\mytag{4.16}\\
\vspace{-1.5ex}
\kern -2em
p&=a_p\,b_p\,c_p,
&q&=a_q\,b_q\,c_q.
\qquad
\endxalignat
$$
The result expressed by the formulas \mythetag{4.15} and 
\mythetag{4.16} is rather important. For this reason it is formulated 
as a lemma.
\mylemma{4.6} For any solution of the Diophantine equations 
\mythetag{3.9}, \mythetag{3.12}, \mythetag{3.13}, and \mythetag{3.14}
with $A_0>0$ and $B_0>0$ there are six positive pairwise mutually 
coprime integer numbers $a_p$, $b_p$, $c_p$, $a_q$, $b_q$, $c_q$ such 
that the numbers $A_0$, $B_0$, $p$, $q$ are expressed through them by 
means of the formulas \mythetag{4.16}.
\endproclaim
     Let's substitute \mythetag{4.16} into the equation \mythetag{4.10}.
As a result we obtain the following equation with respect to the
numbers $a_p$, $b_p$, $c_p$, $a_q$, $b_q$, $c_q$:
$$
\hskip -2em
\gathered
a_p^{10}\,a_q^{10}\,c_p^4\,c_q^4+6\,a_p^6\,a_q^{10}\,c_p^2\,c_q^6\,b_q^4
-a_p^8\,a_q^8\,c_p^4\,c_q^4\,b_q^2\,b_p^2\,-\\
-\,2\,a_p^{10}\,a_q^6\,c_p^6\,c_q^2\,b_p^4
+4\,a_p^6\,a_q^6\,b_p^4\,c_p^4\,b_q^4\,c_q^4
+a_p^{10}\,a_q^2\,b_p^8\,c_p^8\,+\\
+\,a_p^2\,a_q^{10}\,b_q^8\,c_q^8+10\,a_p^4\,a_q^8\,b_p^2\,c_p^2\,b_q^6\,c_q^6
-14\,a_p^8\,a_q^4\,b_p^6\,c_p^6\,b_q^2\,c_q^2=\\
=b_p^{10}\,b_q^{10}\,c_p^4\,c_q^4+6\,b_p^{10}\,b_q^6\,c_p^6\,c_q^2\,a_p^4
-b_p^8\,b_q^8\,c_p^4\,c_q^4\,a_q^2\,a_p^2\,-\\
-2\,b_p^6\,b_q^{10}\,c_p^2\,c_q^6\,a_q^4\,
+4\,b_p^6\,b_q^6\,a_p^4\,c_p^4\,a_q^4\,c_q^4
+b_p^2\,b_q^{10}\,a_q^8\,c_q^8\,+\\
+\,b_p^{10}\,b_q^2\,a_p^8\,c_p^8+10\,b_p^8\,b_q^4\,a_p^6\,c_p^6\,a_q^2\,c_q^2
-14\,b_p^4\,b_q^8\,a_p^2\,c_p^2\,a_q^6\,c_q^6.
\endgathered
\mytag{4.17}
$$
\mylemma{4.7} For a given pair of positive coprime integer numbers $p\neq q$
the system of Diophantine equations \mythetag{3.9}, \mythetag{3.12}, 
\mythetag{3.13}, and \mythetag{3.14} is resolvable if and only if there
are six positive integer numbers $a_p$, $b_p$, $c_p$, $a_q$, $b_q$, $c_q$ 
obeying the equation \mythetag{4.17}, obeying the coprimality conditions
\mythetag{4.15}, and such that $p$ and $q$ are expressed through them
by means of the formulas $p=a_p\,b_p\,c_p$ and $q=a_q\,b_q\,c_q$.
\endproclaim
     The lemma~\mythelemma{4.7} follows from the lemmas~\mythelemma{4.3}
and \mythelemma{4.6} due to the calculations in deriving the equation 
\mythetag{4.17} from the equation \mythetag{4.8}. Combining the 
lemmas~\mythelemma{3.1} and \mythelemma{4.7}, now we derive the following
structural theorem for the solutions of the Diophantine equation
\mythetag{1.2}. 
\mytheorem{4.1} For a given pair of positive coprime integer numbers $p\neq q$
the Diophantine equation \mythetag{1.2} is resolvable with respect to the variable 
$t$ if and only if there are six positive integer numbers $a_p$, $b_p$, $c_p$, 
$a_q$, $b_q$, $c_q$ obeying the equation \mythetag{4.17}, obeying the 
coprimality conditions \mythetag{4.15}, and such that $p$ and $q$ are expressed 
through them by means of the formulas $p=a_p\,b_p\,c_p$ and $q=a_q\,b_q\,c_q$.
Under these conditions the equation \mythetag{1.2} has at least two solutions
given by the formulas
$$
\xalignat 2
&\hskip -2em
t=a_p^2\,c_p\,a_q^2\,c_q,
&&t=-a_p^2\,c_p\,a_q^2\,c_q.
\mytag{4.18}
\endxalignat
$$
\endproclaim
    The structural theorem~\mythetheorem{4.1} is the main result of this paper. 
The formulas \mythetag{4.18} in this theorem are immediate from \mythetag{4.16}. 
\head
5. Some applications of the structural theorem. 
\endhead
    Let's choose $q=1$ and assume that $p$ is some prime number. Then $p$ and $q$
are coprime, i\.\,e\. the relationship \mythetag{4.2} is fulfilled. Applying
the formula $q=a_q\,b_q\,c_q$ from the theorem~\mythetheorem{4.1} to this case, 
we get
$$
\xalignat 3
&\hskip -2em
a_q=1\,
&&b_q=1\,
&&c_q=1.
\mytag{5.1}
\endxalignat
$$
Similarly, applying the formula $p=a_p\,b_p\,c_p$ from the theorem~\mythetheorem{4.1}
and taking into account that $p$ is prime, we derive three options
$$
\xalignat 3
&\hskip -2em
a_p=p,
&&b_p=1,
&&c_p=1;
\mytag{5.2}\\
&\hskip -2em
a_p=1,
&&b_p=p,
&&c_p=1;
\mytag{5.3}\\
&\hskip -2em
a_p=1,
&&b_p=1,
&&c_p=p.
\mytag{5.4}
\endxalignat
$$\par
Substituting \mythetag{5.1} and \mythetag{5.2} into \mythetag{4.17},
we obtain the following equation for $p$:
$$
\hskip -2em
16\,p^2-16\,p^8=0.
\mytag{5.5}
$$
The left hand side of the equation \mythetag{5.5} factorizes as 
$$
\hskip -2em
-16\,p^2\,(p-1)\,(p+1)\,(p^2+p+1)\,(p^2-p+1)=0.
\mytag{5.6}
$$
Therefore the only integer solutions of the equation \mythetag{5.6} 
are
$$
\xalignat 2
&\hskip -2em
p=-1,
&&p=1.
\mytag{5.7}
\endxalignat
$$
The first of them is negative. The second one is positive, but 
$p=1$ contradicts the inequality $p\neq q$ in the 
theorem~\mythetheorem{4.1} since $q=1$ in our 
present case.\par
    The second option is given by the formulas \mythetag{5.3}. 
Substituting \mythetag{5.1} and \mythetag{5.3} into \mythetag{4.17}, 
we obtain another equation for $p$:
$$
\hskip -2em
-8\,p^{10}-8\,p^8-16\,p^6+16\,p^4+8\,p^2+8=0.
\mytag{5.8}
$$
The left hand side of the equation \mythetag{5.8} factorizes as 
follows:
$$
\hskip -2em
-8\,(p-1)\,(p+1)\,(p^8+2\,p^6+4\,p^4+2\,p^2+1)=0. 
\mytag{5.9}
$$
Again the only integer solutions of the equation \mythetag{5.9}
are given by the formulas \mythetag{5.7}. The first of these two
solutions is negative, while the second contradicts the inequality
$p\neq q$ in the theorem~\mythetheorem{4.1}.\par
    And finally, the third option is given by the formulas
\mythetag{5.4}. Substituting \mythetag{5.1} and \mythetag{5.4} into 
\mythetag{4.17}, we obtain the third equation for $p$:
$$
\hskip -2em
-32\,p^6+32\,p^2=0.
\mytag{5.10}
$$
The left hand side of the equation \mythetag{5.10} factorizes as 
follows:
$$
\hskip -2em
-32\,p^2\,(p-1)\,(p+1)\,(p^2+1)=0.
\mytag{5.11}
$$
The integer solutions of the equation \mythetag{5.11} are given by
the formula \mythetag{5.7}. The first of them is negative, while the 
second one contradicts the inequality $p\neq q$. Thus, none of the 
equations \mythetag{5.6}, \mythetag{5.9}, and \mythetag{5.11} has a 
solution suitable for the theorem~\mythetheorem{4.1}. Applying this 
theorem, we can formulate the following result.
\mytheorem{5.1} For $q=1$ and for any positive prime integer $p$ 
the polynomial Diophantine equation \mythetag{1.2} has no integer 
\pagebreak solutions. 
\endproclaim
    Exchanging $p$ and $q$ in the theorem~\mythetheorem{5.1}, we
obtain another theorem.
\mytheorem{5.2} For $p=1$ and for any positive prime integer $q$ 
the polynomial Diophantine equation \mythetag{1.2} has no integer 
solutions. 
\endproclaim
     The proof of the theorem~\mythetheorem{5.2} is similar to the
above proof of the theorem~\mythetheorem{5.1}. The 
theorems~\mythetheorem{5.1} and \mythetheorem{5.2} prove the 
theorem~\mythetheorem{1.1} in two special cases where $p$ is prime
and $q=1$ and where $q$ is prime and $p=1$. Other examples of applying 
the structural theorem~\mythetheorem{4.1} will be given in a separate
paper. 


\Refs
\ref\myrefno{1}\by Sharipov~R.~A.~\paper Perfect cuboids and irreducible 
polynomials\jour e-print \myhref{http://arxiv.org/abs/1108.5348}{arXiv:1108.5348} 
in Electronic Archive \myEarXivlink
\endref
\ref\myrefno{2}\paper
\myhref{http://en.wikipedia.org/wiki/Euler\podcherkivanie 
brick}{Euler brick}\jour Wikipedia\publ 
Wikimedia Foundation Inc.\publaddr San Francisco, USA 
\endref
\ref\myrefno{3}\by Halcke~P.\book Deliciae mathematicae oder mathematisches 
Sinnen-Confect\publ N.~Sauer\publaddr Hamburg, Germany\yr 1719
\endref
\ref\myrefno{4}\by Saunderson~N.\book Elements of algebra, {\rm Vol. 2}\publ
Cambridge Univ\. Press\publaddr Cambridge\yr 1740 
\endref
\ref\myrefno{5}\by Euler~L.\book Vollst\"andige Anleitung zur Algebra
\publ Kayserliche Akademie der Wissenschaften\publaddr St\.~Petersburg
\yr 1771
\endref
\ref\myrefno{6}\by Dickson~L.~E\book History of the theory of numbers, 
{\rm Vol\. 2}: Diophantine analysis\publ Dover\publaddr New York\yr 2005
\endref
\ref\myrefno{7}\by Kraitchik~M.\paper On certain rational cuboids
\jour Scripta Math\.\vol 11\yr 1945\pages 317--326
\endref
\ref\myrefno{8}\by Kraitchik~M.\book Th\'eorie des Nombres,
{\rm Tome 3}, Analyse Diophantine et application aux cuboides 
rationelles \publ Gauthier-Villars\publaddr Paris\yr 1947
\endref
\ref\myrefno{9}\by Kraitchik~M.\paper Sur les cuboides rationelles
\jour Proc\. Int\. Congr\. Math\.\vol 2\yr 1954\publaddr Amsterdam
\pages 33--34
\endref
\ref\myrefno{10}\by Bromhead~T.~B.\paper On square sums of squares
\jour Math\. Gazette\vol 44\issue 349\yr 1960\pages 219--220
\endref
\ref\myrefno{11}\by Lal~M., Blundon~W.~J.\paper Solutions of the 
Diophantine equations $x^2+y^2=l^2$, $y^2+z^2=m^2$, $z^2+x^2
=n^2$\jour Math\. Comp\.\vol 20\yr 1966\pages 144--147
\endref
\ref\myrefno{12}\by Spohn~W.~G.\paper On the integral cuboid\jour Amer\. 
Math\. Monthly\vol 79\issue 1\pages 57-59\yr 1972 
\endref
\ref\myrefno{13}\by Spohn~W.~G.\paper On the derived cuboid\jour Canad\. 
Math\. Bull\.\vol 17\issue 4\pages 575-577\yr 1974
\endref
\ref\myrefno{14}\by Chein~E.~Z.\paper On the derived cuboid of an 
Eulerian triple\jour Canad\. Math\. Bull\.\vol 20\issue 4\yr 1977
\pages 509--510
\endref
\ref\myrefno{15}\by Leech~J.\paper The rational cuboid revisited
\jour Amer\. Math\. Monthly\vol 84\issue 7\pages 518--533\yr 1977
\moreref see also Erratum\jour Amer\. Math\. Monthly\vol 85\page 472
\yr 1978
\endref
\ref\myrefno{16}\by Leech~J.\paper Five tables relating to rational cuboids
\jour Math\. Comp\.\vol 32\yr 1978\pages 657--659
\endref
\ref\myrefno{17}\by Spohn~W.~G.\paper Table of integral cuboids and their 
generators\jour Math\. Comp\.\vol 33\yr 1979\pages 428--429
\endref
\ref\myrefno{18}\by Lagrange~J.\paper Sur le d\'eriv\'e du cuboide 
Eul\'erien\jour Canad\. Math\. Bull\.\vol 22\issue 2\yr 1979\pages 239--241
\endref
\ref\myrefno{19}\by Leech~J.\paper A remark on rational cuboids\jour Canad\. 
Math\. Bull\.\vol 24\issue 3\yr 1981\pages 377--378
\endref
\ref\myrefno{20}\by Korec~I.\paper Nonexistence of small perfect 
rational cuboid\jour Acta Math\. Univ\. Comen\.\vol 42/43\yr 1983
\pages 73--86
\endref
\ref\myrefno{21}\by Korec~I.\paper Nonexistence of small perfect 
rational cuboid II\jour Acta Math\. Univ\. Comen\.\vol 44/45\yr 1984
\pages 39--48
\endref
\ref\myrefno{22}\by Wells~D.~G.\book The Penguin dictionary of curious and 
interesting numbers\publ Penguin publishers\publaddr London\yr 1986
\endref
\ref\myrefno{23}\by Bremner~A., Guy~R.~K.\paper A dozen difficult Diophantine 
dilemmas\jour Amer\. Math\. Monthly\vol 95\issue 1\yr 1988\pages 31--36
\endref
\ref\myrefno{24}\by Bremner~A.\paper The rational cuboid and a quartic surface
\jour Rocky Mountain J\. Math\. \vol 18\issue 1\yr 1988\pages 105--121
\endref
\ref\myrefno{25}\by Colman~W.~J.~A.\paper On certain semiperfect cuboids\jour
Fibonacci Quart.\vol 26\issue 1\yr 1988\pages 54--57\moreref see also\nofrills 
\paper Some observations on the classical cuboid and its parametric solutions
\jour Fibonacci Quart\.\vol 26\issue 4\yr 1988\pages 338--343
\endref
\ref\myrefno{26}\by Korec~I.\paper Lower bounds for perfect rational cuboids 
\jour Math\. Slovaca\vol 42\issue 5\yr 1992\pages 565--582
\endref
\ref\myrefno{27}\by Guy~R.~K.\paper Is there a perfect cuboid? Four squares 
whose sums in pairs are square. Four squares whose differences are square 
\inbook Unsolved Problems in Number Theory, 2nd ed.\pages 173--181\yr 1994
\publ Springer-Verlag\publaddr New York 
\endref
\ref\myrefno{28}\by Rathbun~R.~L., Granlund~T.\paper The integer cuboid table 
with body, edge, and face type of solutions\jour Math\. Comp\.\vol 62\yr 1994
\pages 441--442
\endref
\ref\myrefno{29}\by Van Luijk~R.\book On perfect cuboids, \rm Doctoraalscriptie
\publ Mathematisch Instituut, Universiteit Utrecht\publaddr Utrecht\yr 2000
\endref
\ref\myrefno{30}\by Rathbun~R.~L., Granlund~T.\paper The classical rational 
cuboid table of Maurice Kraitchik\jour Math\. Comp\.\vol 62\yr 1994
\pages 442--443
\endref
\ref\myrefno{31}\by Peterson~B.~E., Jordan~J.~H.\paper Integer hexahedra equivalent 
to perfect boxes\jour Amer\. Math\. Monthly\vol 102\issue 1\yr 1995\pages 41--45
\endref
\ref\myrefno{32}\by Rathbun~R.~L.\paper The rational cuboid table of Maurice 
Kraitchik\jour e-print \myhref{http://arxiv.org/abs/math/0111229}{math.HO/0111229} 
in Electronic Archive \myEarXivlink
\endref
\ref\myrefno{33}\by Hartshorne~R., Van Luijk~R.\paper Non-Euclidean Pythagorean 
triples, a problem of Euler, and rational points on K3 surfaces\publ e-print 
\myhref{http://arxiv.org/abs/math/0606700}{math.NT/0606700} 
in Electronic Archive \myEarXivlink
\endref
\ref\myrefno{34}\by Waldschmidt~M.\paper Open diophantine problems\publ e-print 
\myhref{http://arxiv.org/abs/math/0312440}{math.NT/0312440} 
in Electronic Archive \myEarXivlink
\endref
\ref\myrefno{35}\by Ionascu~E.~J., Luca~F., Stanica~P.\paper Heron triangles 
with two fixed sides\publ e-print \myhref{http://arxiv.org/abs/math/0608185}
{math.NT/0608} \myhref{http://arxiv.org/abs/math/0608185}{185} in Electronic 
Archive \myEarXivlink
\endref
\ref\myrefno{36}\by Sloan~N.~J.~A\paper Sequences 
\myhref{http://oeis.org/A031173}{A031173}, 
\myhref{http://oeis.org/A031174}{A031174}, and \myhref{http://oeis.org/A031175}
{A031175}\jour On-line encyclopedia of integer sequences\publ OEIS Foundation 
Inc.\publaddr Portland, USA
\endref
\ref\myrefno{37}\by Stoll~M., Testa~D.\paper The surface parametrizing cuboids
\jour e-print \myhref{http://arxiv.org/abs/1009.0388}{arXiv:1009.0388} 
in Electronic Archive \myEarXivlink
\endref
\ref\myrefno{38}\by Sharipov~R.~A.~\paper A note on the first cuboid conjecture
\jour e-print \myhref{http://arxiv.org/abs/1109.2534}{arXiv:1109.2534} 
in Electronic Archive \myEarXivlink
\endref
\endRefs
\enddocument
\end